# Divergent Infinite Series – Ramanujan's Initial Intuition

**Mario M. Attard**[1]

**Abstract.** This paper investigates Srinivasa Ramanujan's initial intuitive methodology for assigning the finite value -1/12 to the sum of the divergent infinite series of all positive integers. We systematically examine Ramanujan's initial method, originally sketched in his notebooks, and set the methodology into an algebraic framework. The methodology has limited applicability to other classes of divergent series. The methodology is extended to assign a Ramanujan smoothed sum to the infinite sequences of integers raised to a positive integer power and to figurate binomial number sequences, including triangular numbers, tetrahedral numbers, and higher-dimensional analogues, avoiding analytical continuation. A key finding establishes that the Ramanujan smoothed sums of figurate binomial sequences are intrinsically connected to logarithmic numbers (Gregory coefficients), providing a novel perspective on Ramanujan summation through the lens of classical combinatorial functions. The paper applies asymptotic expansions of associated rational generating functions to demonstrate consistency with established results from analytic continuation methods. The results illuminate the deeper mathematical structures underlying Ramanujan's intuitive insights and suggest new avenues for research in divergent series summation.

# 1 Introduction

Srinivasa Ramanujan was a mathematician known for his exceptional numerical skills and deep intuitive understanding of numerical patterns. He famously assigned the sum of all positive integers to -1/12. Of course, this cannot be the case, as all terms in the summation are positive and the sum is classically divergent with the partial sums growing without bounds. However, the series of all positive integers 1+2+3+4+5+6+…can be associated with a rational generating function applied at the divergent point, for example:

$$\lim_{x\to 1^-}\left[\frac{x^m}{(1-x)^2}=x^m\sum_{n=1}^{\infty}nx^{n-1}=x^m\left(1+2x+3x^2+4x^3+5x^4+6x^5+7x^6+\quad ...\right)\right]\qquad n\in\mathbb{N}\quad x,m\in\mathbb{R}\tag{1}$$

The limit $x\to 1^-$ is outside the radius of convergence of the series, $|x|<1$. Appendix A gives a linear algebra technique using matrix inversion that can be used to derive the generating function for some series.

[1] Honorary Associate Professor, School of Civil and Environmental Engineering, University of New south Wales, Australia. m.attard@unsw.edu.au

Under an asymptotic expansion and/or regulation of Eqn (1), a constant term is revealed at the limit point once the divergent terms are disregarded. For example, substituting the asymptotic gauge function $x=e^{-\lambda}\to 1^{-}\ \left(\lambda\to 0^{+}\right)$ into Eqn (1) reveals:

$$\lim_{x\to 1^{-}}\left[\frac{x^{m}}{\left(1-x\right)^{2}}\right]=\lim_{\lambda\to 0^{+}}\left[\frac{e^{-m\lambda}}{\left(1-e^{-\lambda}\right)^{2}}\right]=\lim_{\lambda\to 0^{+}}\left[\frac{1}{\lambda^{2}}+\frac{1-m}{\lambda}+\left(\frac{5}{12}-m+\frac{m^{2}}{2}\right)+\frac{\lambda}{12}\left(1-5m+6m^{2}-2m^{3}\right)+....\right] \quad (2)$$

If the divergent terms $\left(\frac{1}{\lambda^{2}},\frac{1-m}{\lambda}\right)$ are ignored, as $\left(\lambda\to 0^{+}\right)$ there is the constant term $\left(\frac{5}{12}-m+\frac{m^{2}}{2}\right)$ which can be viewed as "smoothed" values of the sum of the series near the divergence. When *m=1*, the constant term is -1/12 aligning with Ramanujan's result. The series sum is not equal to -1/12 in the conventional sense of summation. The dependence of the constant term $\left(\frac{5}{12}-m+\frac{m^{2}}{2}\right)$ on the index *m* shows the importance of the index shift in the generating function described in Eqn (1) on the resulting asymptotic expansion. Interestingly, the roots of the constant term are $m=1\pm\frac{1}{\sqrt{6}}$ and have zero Ramanujan smoothed sum. For the special case of *m=1*, the Ramanujan smoothed sum can also be represented as:

$$\lim_{x\to 1^{-}}\left[\frac{x}{\left(1-x\right)^{2}}-\frac{1}{\left(\ln\left(x\right)\right)^{2}}\right]=\lim_{x\to 1^{-}}\left[\sum_{n=1}^{\infty}nx^{n}-\frac{1}{\left(\ln\left(x\right)\right)^{2}}\right]=\frac{-1}{12} \quad (3)$$

This is a form of regularizing of an infinite divergent series. One of the most famous examples of regularizing of a divergent series, as pointed out by Kowalenko [1] and unbeknown to Euler at the time, was the equation for the Euler-Mascheroni constant $\gamma$ derived by Euler (see also [2, 3]:

$$\lim_{n\to\infty}\left[1+\frac{1}{2}+\frac{1}{3}+\frac{1}{4}+\frac{1}{5}+...+\frac{1}{n}-\ln\left(n\right)\right]=\gamma \quad (4)$$

Using rigorous methods like zeta regularization, analytic continuation of the Riemann zeta function $\zeta\left(s\right)$, gives the smoothed sum of all positive integers for *s=-1* as $\zeta\left(-1\right)$ = -1/12 matching Ramanujan's initial result.

It is however instructive to review Ramanujan's initial intuitive method where he arrived at a sum of -1/12 and examine whether his initial methodology could be expanded to other sequences. We will first go through Ramanujan's initial method step by step, that was roughly sketched in his notebook [4-7]. Sivaraman [8] also reviewed Ramanujan's initial method obtained from Ramanujan's writing and using elementary concepts avoiding analytical continuation, extended Ramanujan's initial method to the smoothed sum of integers which are the zeta function of negative integers, $\zeta\left(-n\right)\quad n\in\mathbb{Z}^{+}$. In this paper, Ramanujan's initial method is also extended to the sum of integers raised to an integer power but using a different approach to Sivaraman and also avoiding analytical continuation. The methodology examined in this paper is also extended to assign a Ramanujan smoothed sum to infinite divergent Figurate Binomial number sequences (such as the natural numbers, triangular numbers, tetrahedral numbers, pentatope numbers etc.). It is shown that the Ramanujan smooth sum of Figurate

Binomial number sequences are related to Logarithmic numbers (Gregory coefficients). Linking Figurate number sequences to Gregory coefficients provides a new perspective on Ramanujan's summation.

$$\begin{aligned} S_1 &= 1+2+3+4+5+6+7+\ \ldots \overset{\Re}{=} c+\Theta_1 \qquad \Theta_i - \textit{divergent terms} \\ -4S_1 &= \quad -4 \quad -8 \quad -12 \quad \ldots \overset{\Re}{=} -4c-\Theta_2 \\ S_1-4S_1 &= 1-2+3-4+5-6+7+\ \ldots \overset{A}{=} \frac{1}{4} \end{aligned} \tag{5}$$

Ramanujan's initial method begins with writing out the sequence of all positive integers as expressed in the first row of Eqn (5) with summation $S_1$. In this paper, the Ramanujan summation is identified with a modified equal sign "$\overset{\Re}{=}$" with the superscript meaning an assignment based on Ramanujan's summation. The assigned value is "*c*" being a constant. In this paper there is an additional term "$\Theta_i$" placed to identify that there are divergent terms as in the asymptotic expansion of the associated rational generating function in Eqn (2). Next, the sequence of all positive integers is multiplied by -4 and shifted as shown in the second row of Eqn (5). Ramanujan assumes that the smoothed sum of the second row will be -4 times the smoothed sum of the first row. The first two rows are then added together to give Euler's alternating series 1-2+3-4+5-6+…, which can itself be linked to a rational generating function $\lim\limits_{x\to 1^-}\left[\frac{x^m}{(1+x)^2} = x^m\sum\limits_{n=1}^{\infty}(-1)^{n-1}nx^{n-1}\right]$ $n\in\mathbb{N}$. If we substitute an asymptotic gauge function $x=e^{-\lambda}\to 1^-\ (\lambda\to 0^+)$ into this generating function, we see that there are no divergent terms and there is a constant term $\frac{1}{4}$ which corresponds to the Abel summation for Euler's alternating series.

$$\lim_{x\to 1^-}\left[\frac{x^m}{(1+x)^2}\right] = \lim_{\lambda\to 0^+}\left[\frac{e^{-m\lambda}}{\left(1+e^{-\lambda}\right)^2}\right] = \lim_{\lambda\to 0^+}\left[\frac{1}{4}+\frac{\lambda}{4}(1-m)+\frac{\lambda}{16}\left(1-4m+2m^2\right)+\ldots\right]\overset{A}{=}\frac{1}{4} \tag{6}$$

Abel's theorem [1, 9, 10] states that if a power series converges for $|x|<1$ then the sum of the coefficients of the power series at $x=1$ is the limit of the powers series as $x\to 1^-$. The Abel sum of the coefficients of a power series is the limit of a generating function as x approaches 1 from below if the limit exists, therefore $\lim\limits_{x\to 1^-}\left[\frac{x^m}{(1+x)^2}\right]\overset{A}{=}\frac{1}{4}$. The special equal sign $\overset{A}{=}$ is used to identify an Abel summation.

Ramanujan assumes that the first two rows in Eqn (5) are summable. Hence, he obtained (with the divergent terms either disregarded or assumed to cancel each other):

$$c+\Theta_1-4c-\Theta_2 \overset{\Re}{=} \frac{1}{4} \qquad \therefore\ c=-\frac{1}{12} \tag{7}$$

The series in Eqn (5) are recast here as infinite polynomials with integer coefficients being the terms of the series. Possible generating functions are.

$$\frac{x}{(1-x)^2} = x\left(1+2x+3x^2+4x^3+5x^4+6x^5+7x^6+8x^7+\ldots\right) \quad = \sum_{n=1}^{\infty}(n)x^n \qquad |x|<1$$
$$\frac{x^2}{\left(1-x^2\right)^2} = x^2\left(1+2x^2+3x^4+4x^6+5x^8+6x^{10}+7x^{12}+8x^{14}+\ldots\right) \quad = \sum_{n=1}^{\infty}(n)x^{2n} \qquad |x|<1 \tag{8}$$
$$\frac{x}{(1+x)^2} = x\left(1-2x+3x^2-4x^3+5x^4-6x^5+7x^6-8x^7+\ldots\right) \quad = \sum_{n=1}^{\infty}(-1)^{n-1}(n)x^n \qquad |x|<1$$

Taking the limit as $x \to 1^-$ of the generating functions and the equivalent series expansion gives the summation of the integer series. In the second row, the generating function is the first row's generating function but with all *x's* replaced by $x^2$ , while in the third row, the *x* in the denominator is replaced by -*x*. Therefore Eqn (5) could be represented by:

$$\lim_{x\to1^-}\left[x\left(1+2x+3x^2+4x^3+5x^4+6x^5+7x^6+8x^7+\quad\ldots\right)\right] = \lim_{x\to1^-}\left[\frac{x}{(1-x)^2}\right] \overset{\Re}{=} c+\Theta_1$$
$$4\lim_{x\to1^-}\left[x^2\left(\quad -1 \qquad -2x^2 \qquad -3x^4 \qquad -4x^6 \qquad \ldots\right)\right] = \lim_{x\to1^-}\left[\frac{-4x^2}{\left(1-x^2\right)^2}\right] \overset{\Re}{=} -4c-\Theta_2 \tag{9}$$
$$\lim_{x\to1^-}\left[x\left(1-2x+3x^2-4x^3+5x^4-6x^5+7x^6+ \qquad \ldots\right)\right] = \lim_{x\to1^-}\left[\frac{x}{(1+x)^2}\right] \overset{A}{=} \frac{1}{4}$$

The algebraic equivalent of the sum of the three generating functions Eqn (9) is:

$$\frac{x}{(1-x)^2} - \frac{4x^2}{\left(1-x^2\right)^2} = \frac{x}{(1+x)^2} \tag{10}$$

The functions $\frac{x}{(1-x)^2}$ and $\frac{x}{(1+x)^2}$ are reflections about the y axis and have singularities at $x=\pm1$, respectively. The function $\frac{4x^2}{\left(1-x^2\right)^2}$ is zero at x=0 and has singularities at $\pm1$ . A major assumption in Ramanujan's initial method is that the "constant terms" once the divergent terms are ignored are additive for the three series (rows) in Eqn (5) and that the constant term for the two functions $\frac{x}{(1-x)^2}$, $\frac{x^2}{\left(1-x^2\right)^2}$ or first two rows are equivalent.

The asymptotic expansion $x=e^{-\lambda} \to 1^- \quad \left(\lambda \to 0^+\right)$ of the generating functions given in Eqn (11) agrees with this proposition.

$$\lim_{x\to1^-}\left[\frac{x}{(1-x)^2}\right]=\lim_{\lambda\to0^+}\left[\frac{e^{-\lambda}}{\left(1-e^{-\lambda}\right)^2}\right]=\lim_{\lambda\to0^+}\left[\frac{1}{\lambda^2}-\frac{1}{12}+\frac{\lambda^2}{240}-\frac{\lambda^4}{6048}+\frac{\lambda^6}{172800}....\right]\overset{\Re}{=}-\frac{1}{12}$$

$$\lim_{x\to1^-}\left[\frac{-4x^2}{\left(1-x^2\right)^2}\right]=\lim_{\lambda\to0^+}\left[\frac{-4e^{-2\lambda}}{\left(1-e^{-2\lambda}\right)^2}\right]=4\lim_{\lambda\to0^+}\left[-\frac{1}{4\lambda^2}+\frac{1}{12}-\frac{\lambda^2}{60}+\frac{\lambda^4}{378}-\frac{\lambda^6}{2700}....\right]\overset{\Re}{=}\frac{1}{3} \tag{11}$$

$$\lim_{x\to1^-}\left[\frac{x}{(1+x)^2}\right]=\lim_{\lambda\to0^+}\left[\frac{e^{-\lambda}}{\left(1+e^{-\lambda}\right)^2}\right]=\lim_{\lambda\to0^+}\left[\frac{1}{4}-\frac{\lambda^2}{16}+\frac{\lambda^4}{96}-\frac{17\lambda^6}{11520}....\right]\overset{A}{=}\frac{1}{4}$$

# 2 Sum of Integers Raised to an Integer Power $\sum_{n=1}^{\infty}n^k \quad k\in\mathbb{Z}$

Here we will apply Ramanujan’s intuitive method to the summation of integers raised to an integer power. The generating functions for polynomials with coefficients being integers raised to an integer power are given by the following (Carlitz identity [11, 12]):

$$\frac{P_{k-1}(x)}{(1-x)^k}=\sum_{n=1}^{\infty}n^{k-1}x^{n-1} \qquad |x|<1 \qquad k=1,2,3... \tag{12}$$

Note the index “k” sums from 1. In the above, $P_k(x)$ is the *kth* Eulerian Polynomial [11-14] which can be extracted from the following series:

$$\frac{1-x}{1-xe^{t(1-x)}}=\sum_{n=0}^{\infty}P_n(x)\frac{t^n}{n!} \tag{13}$$

The Eulerian Polynomials can also be written using Eulerian numbers $\left\langle {n \atop m} \right\rangle$, Bernoulli polynomials $B_m^{m-n}$ or Gregory polynomials $G_m^{(n)}$ as detailed in Appendix B and stated here as:

$$P_n(x)=\sum_{m=0}^{n-1}\left\langle {n \atop m} \right\rangle x^m=\sum_{m=0}^{n}\frac{n!B_m^{m-n}}{m!}(x-1)^m=\sum_{m=0}^{n}(-1)^m n!(n-m)G_m^{(n)}(x-1)^m \tag{14}$$

A recursive formula which can also be used to generate the Eulerian polynomials is:

$$\mathrm{Li}_1(x)=-\ln(1-x) \quad \mathrm{Li}_{-k}(x)=x\frac{d\mathrm{Li}_{1-k}(x)}{dx} \quad P_k(x)=\frac{\mathrm{Li}_{-k}(x)(1-x)^{(k+1)}}{x} \quad k=0,1,2,3,4.. \tag{15}$$

The $\mathrm{Li}_s(z)=\sum_{n=1}^{\infty}\frac{z^n}{n^s}$ is the polylogarithm function of order s and argument *z*. The first six Eulerian Polynomials are listed below (see [13]):

$$\begin{array}{lcc} P_0(x) & 1 & 1 \\ P_1(x) & 1 & 1 \\ P_2(x) & 1+x & 1+x \\ P_3(x) & 1+4x+x^2 & 1+4x+x^2 \\ P_4(x) & 1+11x+11x^2+x^3 & (1+x)\left(1+10x+x^2\right) \\ P_5(x) & 1+26x+66x^2+26x^3+x^4 & 1+26x+66x^2+26x^3+x^4 \\ P_6(x) & x^5+57x^4+302x^3+302x^2+57x+1 & (1+x)\left(1+56x+246x^2+56x^3+x^4\right) \end{array} \qquad (16)$$

Note the following well known identities for Eulerian polynomials:

$$P_{k-1}(1)=(k-1)! \qquad P_{k-1}(-1)=\frac{(-1)^{(k-1)}\,2^k\left(2^k-1\right)B_k}{k} \qquad k=1,2,3... \qquad (17)$$

In which $B_k$ is the *k*th Bernoulli number defined by the series expansion of $\frac{t}{e^t-1}=\sum_{n=0}^{\infty}\frac{B_n}{n!}t^n \quad |t|<2\pi$.

Appendix B gives general formula for higher order derivatives of Eulerian Polynomials evaluated at *x=1*.

Before looking at the general case for the summation of integers raised to an integer power, let us consider the summation of the integers cubed. For this we need $P_3(x)$ and the generating function given below.

$$\frac{x\left(1+4x+x^2\right)}{(1-x)^4}=\sum_{n=1}^{\infty}(n)^3x^n \qquad |x|<1 \qquad (18)$$

The twisted or alternating counterpart generating function is:

$$\frac{x\left(1-4x+x^2\right)}{(1+x)^4}=\sum_{n=1}^{\infty}(-1)^n(n)^3x^n \qquad |x|<1 \qquad (19)$$

The difference between the two generating functions is:

$$\frac{x\left(1+4x+x^2\right)}{(1-x)^4}-\frac{x\left(1-4x+x^2\right)}{(1+x)^4}=\frac{16x^2\left(1+4x^2+x^4\right)}{\left(1-x^2\right)^4}=16\sum_{n=1}^{\infty}(n)^3x^{2n} \qquad |x|<1 \qquad (20)$$

We note that the resulting squared generating function is multiplied by a constant, 16. Therefore, we can write using Ramanujan's methodology and Abel summation:

$$\lim_{x\to 1^-}\left[\frac{x\left(1+4x+x^2\right)}{\left(1-x\right)^4}\right]=\sum_{n=1}^{\infty}\left(n\right)^3 \overset{\Re}{=} c+\Theta_1 \qquad \lim_{x\to 1^-}\left[\frac{x\left(1-4x+x^2\right)}{\left(1+x\right)^4}\right]=\sum_{n=1}^{\infty}\left(-1\right)^n\left(n\right)^3 \overset{A}{=} \frac{-1}{8}$$

$$\lim_{x\to 1^-}\left[\frac{16x^2\left[1+4x^2+x^4\right]}{\left(1-x^2\right)^4}\right]=\sum_{n=1}^{\infty}16\left(n\right)^3 \overset{\Re}{=} 16c+\Theta_2 \tag{21}$$

Therefore, the Ramanujan summation is $\frac{1}{120}$,

$$c+\Theta_1-16c-\Theta_2 \overset{\Re}{=} \frac{-1}{8} \qquad \therefore\ c=\frac{1}{120} \qquad \therefore\ \sum_{n=1}^{\infty}n^3 \overset{\Re}{=} \frac{1}{120} \tag{22}$$

The asymptotic expansion $x=e^{-\lambda}\to 1^-\left(\lambda\to 0^+\right)$ of the generating functions given in Eqn (21) agrees with this proposition.

$$\lim_{x\to 1^-}\left[\frac{x\left(1+4x+x^2\right)}{\left(1-x\right)^4}\right]=\lim_{\lambda\to 0^+}\left[\frac{6}{\lambda^4}+\frac{1}{120}-\frac{\lambda^2}{504}+....\right] \overset{\Re}{=} \frac{1}{120}$$

$$\lim_{x\to 1^-}\left[\frac{16x^2\left[1+4x^2+x^4\right]}{\left(1-x^2\right)^4}\right]=\lim_{\lambda\to 0^+}\left[\frac{6}{\lambda^4}+\frac{16}{120}-\frac{8\lambda^2}{63}+.......\right] \overset{\Re}{=} \frac{16}{120} \tag{23}$$

$$\lim_{x\to 1^-}\left[\frac{x\left(1-4x+x^2\right)}{\left(1+x\right)^4}\right]=\lim_{\lambda\to 0^+}\left[-\frac{1}{8}+\frac{\lambda^2}{8}-\frac{17\lambda^4}{384}+....\right] \overset{A}{=} -\frac{1}{8}$$

From the above we can observe that the regularization of the divergent series is:

$$\lim_{x\to 1^-}\left[\frac{x\left(1+4x+x^2\right)}{\left(1-x\right)^4}-\frac{6}{\left(-\ln\left(x\right)\right)^4}\right]=\lim_{x\to 1^-}\left[\sum_{n=1}^{\infty}\left(n\right)^3x^n-\frac{6}{\left(-\ln\left(x\right)\right)^4}\right]=\frac{1}{120} \tag{24}$$

For the general case $\sum_{n=1}^{\infty}n^k \quad k\in\mathbb{Z}$, the generating function and the twisted counterpart are:

$$\frac{xP_{k-1}\left(x\right)}{\left(1-x\right)^k}=\sum_{n=1}^{\infty}n^{k-1}x^n \qquad \frac{xP_{k-1}\left(-x\right)}{\left(1+x\right)^k}=\sum_{n=1}^{\infty}\left(-1\right)^{n+1}n^{k-1}x^n \qquad k=1,2,3... \tag{25}$$

The difference is then:

$$\frac{\left[xP_{k-1}\left(x\right)\left(1+x\right)^k-xP_{k-1}\left(-x\right)\left(1-x\right)^k\right]}{\left(1-x\right)^k\left(1+x\right)^k}=\frac{2^kx^2P_{k-1}\left(x^2\right)}{\left(1-x^2\right)^k}=\sum_{n=1}^{\infty}2^kn^{k-1}x^{2n} \qquad k=1,2,3... \tag{26}$$

The odd power terms of $xP_{k-1}(x)(1+x)^k$ are cancelled by the odd power terms in $-xP_{k-1}(-x)(1-x)^k$ leaving only even powers with a common factor of $x^2$. The resulting polynomial is $2^k$ times the squared generating function. Therefore, using Eqn (17) we can assert:

$$\lim_{x\to 1^-}\left[\frac{xP_{k-1}(x)}{(1-x)^k}\right]=\sum_{n=1}^{\infty} n^{k-1} \overset{\Re}{=} c+\Theta_1 \qquad \lim_{x\to 1^-}\left[\frac{xP_{k-1}(-x)}{(1+x)^k}\right]=\sum_{n=1}^{\infty}(-1)^{n+1} n^{k-1} \overset{A}{=} \frac{P_{k-1}(-1)}{2^k}$$

$$\lim_{x\to 1^-}\left[\frac{2^k x^2 P_{k-1}\left(x^2\right)}{\left(1-x^2\right)^k}\right]=\sum_{n=1}^{\infty} 2^k n^{k-1} \overset{\Re}{=} 2^k c+\Theta_2 \tag{27}$$

$$\therefore \quad c+\Theta_1-2^k c-\Theta_2 \overset{\Re}{=} \frac{P_{k-1}(-1)}{2^k} \qquad \therefore \quad c=\frac{P_{k-1}(-1)}{2^k\left(1-2^k\right)}=\frac{(-1)^{(k-1)}B_k}{k} \qquad k=1,2,3..$$

Hence:

$$\therefore \sum_{n=1}^{\infty} n^{k-1} \overset{\Re}{=} \frac{(-1)^{(k-1)}B_k}{k} \qquad \sum_{n=1}^{\infty}(-1)^{n+1} n^{k-1} \overset{A}{=} \frac{(-1)^{(k-1)}B_k}{k}\left(1-2^k\right) \qquad k=1,2,3.. \tag{28}$$

Table 1 gives some values for the sum of integers raised to an integer power obtained using the extend Ramanujan intuitive method Eqns (28) which agrees with known values for $\zeta(1-k)$, $k$ being a positive integer. The regularization of the divergent series is:

$$\lim_{x\to 1^-}\left[\frac{xP_{k-1}(x)}{(1-x)^k}-\frac{(k-1)!}{(-\ln(x))^k}\right]=\lim_{x\to 1^-}\left[\sum_{n=1}^{\infty} n^{k-1}x^n-\frac{(k-1)!}{(-\ln(x))^k}\right]=\frac{(-1)^{(k-1)}B_k}{k} \qquad k=1,2,3... \tag{29}$$

To prove this limit substitute $1-x=\varepsilon$ in the equation above and the look at the series expansion of the two functions. Using Gregory polynomials of $G_m^{(k)}$ defined in Eqn (90) to (94) in Appendix C, we can write for the $-\dfrac{(k-1)!}{(-\ln(1-\varepsilon))^k}$ expansion:

$$\lim_{\varepsilon\to 0}\left[-\frac{(k-1)!}{(-\ln(1-\varepsilon))^k}\right]=\lim_{\varepsilon\to 0}\left[-k!\,\varepsilon^{-k}\left[\frac{\left(\dfrac{-\ln(1-\varepsilon)}{\varepsilon}\right)^{-k}}{k}\right]\right]=\lim_{\varepsilon\to 0}\left[-k!\sum_{m=0}^{\infty} G_m^{(k)}\varepsilon^{m-k}\right]$$

$$\begin{aligned}
&=\lim_{\varepsilon\to 0}\left[-(k-1)!\,\varepsilon^{-k}-k!\sum_{m=1}^{k-1} G_m^{(k)}\varepsilon^{m-k}-k!\,G_k^{(k)}\varepsilon^{0}-k!\sum_{m=k+1}^{\infty} G_m^{(k)}\varepsilon^{m-k}\right]\\
&=\lim_{\varepsilon\to 0}\left[-(k-1)!\,\varepsilon^{-k}-k!\sum_{m=1}^{k-1} G_m^{(k)}\varepsilon^{m-k}-k!\,G_k^{(k)}\varepsilon^{0}\right] \qquad k=1,2,3...
\end{aligned} \tag{30}$$

For the other function in Eqn (29) $\frac{xP_{k-1}(x)}{(1-x)^k}$, it is convenient to write out the Eulerian polynomials as a power series expansion in terms of derivatives and the function evaluated at the expansion point. Incorporating the identity in Eqn (17) gives:

$$\begin{aligned}
\lim_{x\to1^-}\left[\frac{xP_{k-1}(x)}{(1-x)^k}\right] &= \lim_{\varepsilon\to0}\left[(1-\varepsilon)P_{k-1}(1-\varepsilon)\varepsilon^{-k}\right] \qquad k=1,2,3...\\
&= \lim_{\varepsilon\to0}\left[\varepsilon^{-k}\,P_{k-1}(1-\varepsilon)\Big|_{\varepsilon=0}+\sum_{m=1}^{k-1}\frac{\varepsilon^{m-k}}{m!}\left[\frac{d^mP_{k-1}(1-\varepsilon)}{d\varepsilon^m}-m\frac{d^{m-1}P_{k-1}(1-\varepsilon)}{d\varepsilon^{m-1}}\right]_{\varepsilon=0}\right]\\
&= \lim_{\varepsilon\to0}\left[(k-1)!\varepsilon^{-k}+\sum_{m=1}^{k-1}\frac{\varepsilon^{m-k}}{m!}\left[\frac{d^mP_{k-1}(1-\varepsilon)}{d\varepsilon^m}-m\frac{d^{m-1}P_{k-1}(1-\varepsilon)}{d\varepsilon^{m-1}}\right]_{\varepsilon=0}\right]
\end{aligned} \tag{31}$$

The terms involving derivatives can be represented by Gregory polynomials using equations (88) and (92) therefore:

$$\sum_{m=1}^{k-1}\frac{\varepsilon^{m-k}}{m!}\left[\frac{d^mP_{k-1}(1-\varepsilon)}{d\varepsilon^m}-m\frac{d^{m-1}P_{k-1}(1-\varepsilon)}{d\varepsilon^{m-1}}\right]_{\varepsilon=0}=k!\sum_{m=1}^{k-1}G_m^{(k)}\varepsilon^{m-k} \qquad k=1,2,3... \tag{32}$$

Hence

$$\lim_{x\to1^-}\left[\frac{xP_{k-1}(x)}{(1-x)^k}\right]=\lim_{\varepsilon\to0}\left[(k-1)!\varepsilon^{-k}+k!\sum_{m=1}^{k-1}G_m^{(k)}\varepsilon^{m-k}\right] \qquad k=1,2,3... \tag{33}$$

Combining Eqns (30) and (33), and using Eqn (97) gives the following which confirms Eqn (29):

$$\lim_{x\to1^-}\left[\frac{xP_{k-1}(x)}{(1-x)^k}-\frac{(k-1)!}{(-\ln(x))^k}\right]=-k!G_k^{(k)}=\frac{(-1)^{(k-1)}B_k}{k} \qquad k=1,2,3... \tag{34}$$

Equation (29) agrees with the general regularization of the Riemann Zeta function. Equation 25.12.12 in [15] is:

$$\mathrm{Li}_s(x)-\frac{\Gamma(1-s)}{(-\ln(x))^{1-s}}=\sum_{n=0}^{\infty}\zeta(s-n)\frac{(\ln(x))^n}{n!} \qquad s\in\mathbb{C},\quad s\neq1,2,3.....\ |\ln(x)|<2\pi \tag{35}$$

Where $\Gamma(s)$ is the Gamma function. Taking the asymptotic expansion of Eqn (35) gives the general regularization of the Riemann Zeta function:

$$\begin{aligned}
&\lim_{x\to1^-}\left[\sum_{n=1}^{\infty}n^{-s}x^n-\frac{\Gamma(1-s)}{(-\ln(x))^{1-s}}\right]=\lim_{x\to1^-}\left[\mathrm{Li}_s(x)-\frac{\Gamma(1-s)}{(-\ln(x))^{1-s}}\right]=\zeta(s) \quad \Re s<1\\
&\mathrm{Li}_s(1)=\zeta(s) \qquad \Re s>1
\end{aligned} \tag{36}$$

Table 1 Sum of Integers Raised to an Integer Power

| *k-1* | Ramanujan Smoothed Sum $\frac{-B(k)}{k}$ | $\sum_{n=1}^{\infty} n^{k-1}$ | Generating Function $\frac{xP_{k-1}(x)}{(1-x)^k}$ |
|---|---|---|---|
| *0* | *-1/2* | 1 1 1 1 1 1 1…. | $\frac{x}{(1-x)}$ |
| *1* | *-1/12* | 1 2 3 4 5 6 7…. | $\frac{x}{(1-x)^2}$ |
| *2* | *0* | 1 4 9 16 25 36 49…. | $\frac{x(1+x)}{(1-x)^3}$ |
| *3* | *1/120* | 1 8 27 64 125 216 343…. | $\frac{x(1+4x+x^2)}{(1-x)^4}$ |
| *4* | *0* | 1 16 81 256 625 1296 …. | $\frac{x(1+x)(1+10x+x^2)}{(1-x)^5}$ |
| *5* | *-1/252* | 1 32 243 1024 3125 7776 …. | $\frac{x(1+26x+66x^2+26x^3+x^4)}{(1-x)^6}$ |
| *6* | *0* | 1 64 729 4096 15625 …. | $\frac{x(1+x)(1+56x+246x^2+56x^3+x^4)}{(1-x)^7}$ |
| *7* | *1/240* | 1 128 2187 16384 78125 …. | $\frac{x(1+120x+1191x^2+2416x^3+1191x^4+120x^5+1)}{(1-x)^8}$ |

# 3 Figurate Binomial Sequences

Figurate binomial numbers [16] are an important class of combinational sequences associated with geometric interpretations. Figurate binomial number sequences are defined by the binomial expression:

$$\sum_{n=0}^{\infty}\binom{n+k-1}{n}=\sum_{n=0}^{\infty}\frac{(n+k-1)!}{n!(k-1)!}=\frac{1}{(k-1)!}\prod_{i=0}^{k-2}(n+i) \tag{37}$$

where the binomial coefficients are weighted by figurate numbers. Figurate numbers appear in combinatorial geometry, probability and algebraic combinatorics [16]. Table 2 lists the first seven examples of Figurate number sequences.

The generating function for Figurate number sequences selected here, is $\frac{x}{(1-x)^k}$ while the twisted counterpart is $\frac{x}{(1+x)^k}$. Table 3 gives the difference of the generating function and it's twisted counterpart as while as it's decomposition. The Ramanujan smooth sum for *k=1 & 2* are the same as for the $\sum_{n=1}^{\infty} n^0$ & $\sum_{n=1}^{\infty} n$ which were given previously as $-\frac{1}{2}$ & $-\frac{1}{12}$ , respectively.

Let us firstly consider the Ramanujan summation of the Triangular integers, as an example of the application of Ramanujan's intuitive method to these types of sequences. The generating function and the twisted or alternating counterpart are:

$$\frac{x}{(1-x)^3} = \sum_{n=1}^{\infty} \frac{n(n+1)}{2} x^n \qquad \frac{x}{(1+x)^3} = \sum_{n=1}^{\infty} (-1)^{n-1} \frac{n(n+1)}{2} x^n \qquad |x| < 1 \tag{38}$$

The difference between the generating function and the twisted counterpart is:

$$\frac{x}{(1-x)^3} - \frac{x}{(1+x)^3} = \frac{2x^2\left(x^2+3\right)}{\left(1-x^2\right)^3} = \frac{8x^2}{\left(1-x^2\right)^3} - \frac{2x^2}{\left(1-x^2\right)^2} \tag{39}$$

The squared generating function $\frac{x^2}{\left(1-x^2\right)^3}$ is multiplied by a polynomial rather than a constant. This polynomial is decomposed into partial fractions as in Eqn (39). The Ramanujan smoother sum for the second term in Eqn (39), $\frac{2x^2}{\left(1-x^2\right)^2}$ is known from the smoothed sum calculation of the positive integers.

Table 2 Figurate Binomial Sequence

| *k* | $\sum_{n=0}^{\infty}\binom{n+k-1}{n}$ | | $\frac{x}{(1-x)^k}$ |
|---|---|---|---|
| *1* | 1 1 1 1 1 1 1…. | Unit Integers | $\frac{x}{(1-x)}$ |
| *2* | 1 2 3 4 5 6 7…. | Positive Integers | $\frac{x}{(1-x)^2}$ |
| *3* | 1 3 6 10 15 21 28…. | Triangular Integers | $\frac{x}{(1-x)^3}$ |
| *4* | 1 4 10 20 35 56 84…. | Tetrahedral Integers | $\frac{x}{(1-x)^4}$ |
| *5* | 1 5 15 35 70 126 210…. | Pentatope Integers | $\frac{x}{(1-x)^5}$ |
| *6* | 1 6 21 56 126 252 462…. | Hexateron Integers | $\frac{x}{(1-x)^6}$ |
| *7* | 1 7 28 84 210 462 924…. | Heptapeton Integers | $\frac{x}{(1-x)^7}$ |

Therefore, we can write:

$$\lim_{x\to1^-}\left[\frac{x}{(1-x)^3}\right]=\sum_{n=1}^{\infty}\frac{(n)(n+1)}{2}\overset{\Re}{=}c+\Theta_1 \qquad \lim_{x\to1^-}\left[\frac{x}{(1+x)^3}\right]=\sum_{n=1}^{\infty}(-1)^{n-1}\frac{(n)(n+1)}{2}\overset{A}{=}\frac{1}{8}$$

$$\lim_{x\to1^-}\left[\frac{2x^2\left(x^2+3\right)}{\left(1-x^2\right)^3}\right]=\lim_{x\to1^-}\left[\frac{8x^2}{\left(1-x^2\right)^3}-\frac{2x^2}{\left(1-x^2\right)^2}\right]\overset{\Re}{=}8c+\frac{2}{12}+\Theta_2 \tag{40}$$

Hence,

$$c+\Theta_1-8c-\frac{1}{6}-\Theta_2\overset{\Re}{=}\frac{1}{8}\quad\rightarrow\quad c=\frac{-1}{24}$$

$$\therefore$$

$$\sum_{n=1}^{\infty}\frac{(n)(n+1)}{2}\overset{\Re}{=}\frac{-1}{24}\qquad\lim_{x\to1^-}\left[\frac{x}{(1-x)^3}\right]\overset{\Re}{=}\frac{-1}{24}\qquad\lim_{x\to1^-}\left[\frac{x^2}{\left(1-x^2\right)^3}\right]\overset{\Re}{=}\frac{-1}{24} \tag{41}$$

Ramanujan's smoothed sum can be calculated in a progressive manner for higher order values of *k*, each depending on the other previous calculations. Table 3 gives examples of the calculations up to $k=7$. Table 4 summaries the Ramanujan's smoothed sum which aligns with the asymptotic expansion method. Also listed are the first seven Gregory coefficients (or logarithmic numbers) given by:

$$G_n = 1, \frac{1}{2}, -\frac{1}{12}, \frac{1}{24}, -\frac{19}{720}, \frac{3}{160}, -\frac{863}{60480}, \ldots \qquad n = 0,1,2,3.... \tag{42}$$

The denominator of the Gregory coefficients are related to the Hirzebruch numbers $h_k$ [17, 18]:

$$h_k = 1,\ 2,\ 12,\ 24,\ 720,\ 1440,\ 60480, .... = \prod_{prime\ p \le k+1} p^{\left\lfloor \frac{k}{p-1} \right\rfloor} \tag{43}$$

The Gregory coefficients are also the coefficients of the following generating functions and are related to the multiplicative inverse of the Harmonic series and Bernoulli numbers:

$$\begin{aligned}
&\frac{-x}{\ln(1-x)} = \sum_{n=0}^{\infty} (-1)^n G_n x^n && \frac{x}{\ln(1+x)} = \sum_{n=0}^{\infty} G_n x^n \qquad |x|<1 \\
&\frac{\ln(1-x)}{-x} = \sum_{n=0}^{\infty} \frac{1}{n+1} x^n \quad |x|<1 && t = \ln(1-x) \quad \to \frac{t}{e^t - 1} = \sum_{n=0}^{\infty} \frac{B_n}{n!} t^n \qquad |t| < 2\pi \\
&\frac{\ln(1+x)}{x} = \sum_{n=0}^{\infty} (-1)^n \frac{1}{n+1} x^n \quad |x|<1 && t = \ln(1+x) \quad \to \frac{t}{e^t - 1} = \sum_{n=0}^{\infty} \frac{B_n}{n!} t^n \qquad |t| < 2\pi
\end{aligned} \tag{44}$$

A recursive formula can be derived by looking at the convolution relations between coefficients in the series expansion of $\frac{-x}{\ln(1-x)}$ and its reciprocal $\frac{\ln(1-x)}{-x}$, therefore:

$$\bar{G}_n = -\sum_{k=1}^{n} \frac{\bar{G}_{n-k}}{(k+1)} \quad n>0 \ \ with\ \bar{G}_0 = 1 \qquad G_0 = 1 \quad G_n = (-1)^n \bar{G}_n \quad n = 1,2,3... \tag{45}$$

The Gregory coefficients are also referred to as Bernoulli numbers of the second kind [19], reciprocal logarithmic numbers or Cauchy numbers of the first kind. It is shown here that the Ramanujan summation of the Figurate Binomial Sequences is related to signed Gregory coefficients.

$$\lim_{x \to 1^-} \left[ \frac{x}{(1-x)^k} \right] = \lim_{x \to 1^-} \left[ \sum_{n=1}^{\infty} \binom{n+k-1}{n} x^n \right] \overset{\Re}{=} \bar{G}_k = (-1)^k G_k \qquad k = 1,2,3.. \tag{46}$$

To prove this, substitute the asymptotic expansion $x = e^{-\lambda} \to 1^- \left(\lambda \to 0^+\right)$ into the generating function, that is:

$$\lim_{x \to 1^-} \left[ \frac{x}{(1-x)^k} \right] = \lim_{\lambda \to 0^+} \left[ \frac{e^{-\lambda}}{\left(1-e^{-\lambda}\right)^k} \right] = \lim_{\lambda \to 0^+} \left[ \left( \frac{\lambda}{e^{\lambda} - 1} \right)^k e^{\lambda(k-1)} \lambda^{-k} \right] = \lim_{\lambda \to 0^+} \left[ \sum_{n=0}^{\infty} \frac{B_n^{(k)}(k-1)}{n!} \lambda^{n-k} \right] \tag{47}$$

In the above, $B_n^{(k)}(k-1)$ are the extended generalized Bernoulli polynomials of order k, [15, 20] in the variable $(x = k-1)$. The extended generalized Bernoulli polynomials of order *k*, $B_n^{(k)}(x)$ and the generalized Bernoulli numbers of order *k*, $B_n^{(k)}$ are defined by the following series expansion [15]:

$$\left(\frac{t}{e^t-1}\right)^k = \sum_{n=0}^{\infty}\frac{B_n^{(k)}}{n!}t^n \qquad \left(\frac{t}{e^t-1}\right)^k e^{xt} = \sum_{n=0}^{\infty}\frac{B_n^{(k)}(x)}{n!}t^n \qquad |t|<2\pi \tag{48}$$

By setting $t$=-$t$ it is straight forward to show that:

$$\frac{B_n^{(k)}(k)}{n!} = (-1)^n\frac{B_n^{(k)}}{n!} \tag{49}$$

To determine the relationship between $B_n^{(k)}(k-1)$ and the Gregory coefficients, we begin by establishing the following identity:

$$(k-1)B_n^{(k)}(k-1) = (k-n-1)B_n^{(k-1)}(k-1) \tag{50}$$

This can be written in summation form as:

$$\begin{aligned}
(k-1)\sum_{n=0}^{\infty}B_n^{(k)}(k-1)\frac{t^n}{n!} &= (k-1)\sum_{n=0}^{\infty}B_n^{(k-1)}(k-1)\frac{t^n}{n!} - \sum_{n=0}^{\infty}nB_n^{(k-1)}(k-1)\frac{t^n}{n!} \\
(k-1)\sum_{n=0}^{\infty}B_n^{(k)}(k-1)\frac{t^n}{n!} &= (k-1)\sum_{n=0}^{\infty}B_n^{(k-1)}(k-1)\frac{t^n}{n!} - t\frac{d}{dt}\left(\sum_{n=0}^{\infty}B_n^{(k-1)}(k-1)\frac{t^n}{n!}\right)
\end{aligned} \tag{51}$$

Note the derivative identity:

$$\frac{d}{dt}\left(\sum_{n=0}^{\infty}\frac{B_n^{(k-1)}(k-1)}{n!}t^n\right) = \left(\frac{k-1}{t}\right)\left(1-\frac{t}{e^t-1}\right)\sum_{n=0}^{\infty}\frac{B_n^{(k-1)}(k-1)}{n!}t^n \tag{52}$$

Substituting the above equation into Eqn (51) proves Eqn (50):

$$\begin{aligned}
&(k-1)\sum_{n=0}^{\infty}B_n^{(k)}(k-1)\frac{t^n}{n!} = (k-1)\sum_{n=0}^{\infty}B_n^{(k-1)}(k-1)\frac{t^n}{n!} - (k-1)\sum_{n=0}^{\infty}\frac{B_n^{(k-1)}(k-1)}{n!}t^n \\
&+(k-1)\left(\frac{t}{e^t-1}\right)\sum_{n=0}^{\infty}\frac{B_n^{(k-1)}(k-1)}{n!}t^n \\
&\therefore\quad (k-1)\sum_{n=0}^{\infty}B_n^{(k)}(k-1)\frac{t^n}{n!} = (k-1)\sum_{n=0}^{\infty}B_n^{(k)}(k-1)\frac{t^n}{n!}
\end{aligned} \tag{53}$$

Substituting $n$=$k$ into Eqn (50) and using Eqns (49), we have:

$$(k-1)\frac{B_k^{(k)}(k-1)}{k!} = -\frac{B_k^{(k-1)}(k-1)}{k!} = (-1)^{k+1}\frac{B_k^{(k-1)}}{k!} \tag{54}$$

The *kth* Ramanujan's smoothed sum can be extracted by obtaining the $n^{th}$ term in Eqn (47) and setting $n$=$k$, and incorporating Eqns (54) and (98), resulting in:

$$at \quad n=k \qquad \frac{B_k^{(k)}(k-1)}{k!}=(-1)^k G_k \qquad k=1,2,3... \tag{55}$$

Confirming Eqn (46).

The generating function for Figurate Binomial Sequences is generalised introducing an index shift parameter "*m*" as:

$$\frac{x^m}{(1-x)^k}=\sum_{n=1}^{\infty}\binom{n+k-1}{n}x^{n+m-1} \qquad |x|<1 \qquad n,k\in\mathbb{N} \tag{56}$$

Substituting the asymptotic expansion $x=e^{-\lambda}\to 1^-\left(\lambda\to 0^+\right)$ into Eqn (56) gives:

$$\lim_{x\to 1^-}\left[\frac{x^m}{(1-x)^k}\right]=\lim_{\lambda\to 0^+}\left[\frac{e^{-\lambda m}}{\left(1-e^{-\lambda}\right)^k}\right]=\lim_{\lambda\to 0^+}\left[\frac{\lambda^k}{\left(e^{\lambda}-1\right)^k}e^{\lambda(k-m)}\lambda^{-k}\right]=\lim_{\lambda\to 0^+}\left[\sum_{n=0}^{\infty}\frac{B_n^{(k)}(k-m)}{n!}\lambda^{n-k}\right] \tag{57}$$

Extracting the constant term at the limit gives polynomials in *m* for the Ramanujan summation:

$$\frac{B_k^{(k)}(k-m)}{k!} \tag{58}$$

From the reflection law given below,

$$B_k^{(k)}(k-m)=(-1)^k B_k^{(k)}(m) \tag{59}$$

When *k* is odd, one of the solutions for zero Ramanujan summation $B_k^{(k)}(k-m)=0$ is $m=k/2$ for the shift parameter.

Table 3 Ramanujan's intuitive Method Calculations

| *k* | $\frac{x}{(1-x)^k}-\frac{x}{(1+x)^k}$ | Decomposition | Ramanujan Smoothed Sum Calculation |
|---|---|---|---|
| *1* | $\frac{2x^2}{\left(1-x^2\right)}$ | | *-1/2* |
| *2* | $\frac{4x^2}{\left(1-x^2\right)^2}$ | | *-1/12* |
| *3* | $\frac{2x^2\left(x^2+3\right)}{\left(1-x^2\right)^3}$ | $\frac{8x^2}{\left(1-x^2\right)^3}-\frac{2x^2}{\left(1-x^2\right)^2}$ | $c+\Theta_1-8c-\frac{1}{6}-\Theta_2 \overset{\Re}{=} \frac{1}{8}$ <br> $\rightarrow \quad c=\frac{-1}{24}$ |
| *4* | $\frac{8x^2\left(x^2+1\right)}{\left(1-x^2\right)^4}$ | $\frac{16x^2}{\left(1-x^2\right)^4}-\frac{8x^2}{\left(1-x^2\right)^3}$ | $c+\Theta_1-16c-\frac{8}{24}-\Theta_2 \overset{\Re}{=} \frac{1}{16}$ <br> $\rightarrow \quad c=\frac{-19}{720}$ |
| *5* | $\frac{2x^2\left(x^4+10x^2+5\right)}{\left(1-x^2\right)^5}$ | $\frac{32x^2}{\left(1-x^2\right)^5}-\frac{24x^2}{\left(1-x^2\right)^4}+\frac{2x^2}{\left(1-x^2\right)^3}$ | $c+\Theta_1-31c-24\left(\frac{-19}{720}\right)-2\left(\frac{-1}{24}\right)+\Theta_2$ <br> $\overset{\Re}{=} \frac{1}{32} \quad \rightarrow \quad c=\frac{-27}{1440}$ |
| *6* | $\frac{4x^2\left(x^2+3\right)\left(3x^2+1\right)}{\left(1-x^2\right)^6}$ | $\frac{64x^2}{\left(1-x^2\right)^6}-\frac{64x^2}{\left(1-x^2\right)^5}+\frac{12x^2}{\left(1-x^2\right)^4}$ | $c+\Theta_1-64c+64\left(\frac{-3}{160}\right)-12\left(\frac{-19}{720}\right)+\Theta_2$ <br> $\overset{\Re}{=} \frac{1}{64} \quad \rightarrow \quad c=\frac{-863}{60480}$ |
| *7* | $\frac{2x^2\left(x^6+21x^4+35x^2+7\right)}{\left(1-x^2\right)^7}$ | $\frac{128x^2}{\left(1-x^2\right)^7}-\frac{160x^2}{\left(1-x^2\right)^6}+\frac{48x^2}{\left(1-x^2\right)^5}$ <br> $-\frac{2x^2}{\left(1-x^2\right)^4}$ | $c+\Theta_1-128c+160\left(\frac{-863}{60480}\right)-48\left(\frac{-3}{160}\right)$ <br> $+2\left(\frac{-19}{720}\right)+\Theta_2 \overset{\Re}{=} \frac{1}{128} \quad \rightarrow \quad c=\frac{-1375}{120960}$ |

Table 4 Ramanujan's Intuitive Method Summation for Figurate Binomial Sequences

| *k* | Ramanujan Smoothed Sum | Figurate Binomial Sequence $\sum_{n=0}^{\infty}\binom{n+k-1}{n}$ | $\frac{x}{(1-x)^k}$ | Asymptotic Expansion | Gregory Coefficients $(-1)^k G_k$ |
|---|---|---|---|---|---|
| *1* | *-1/2* | 1 1 1 1 1 1 1…. | $\frac{1}{(1-x)}$ | *-1/2* | *-1/2* |
| *2* | *-1/12* | 1 2 3 4 5 6 7…. | $\frac{x}{(1-x)^2}$ | *-1/12* | *-1/12* |
| *3* | *-1/24* | 1 3 6 10 15 21 28…. | $\frac{x}{(1-x)^3}$ | *-1/24* | *-1/24* |
| *4* | *-19/720* | 1 4 10 20 35 56 84…. | $\frac{x}{(1-x)^4}$ | *-19/720* | *-19/720* |
| *5* | *-27/1440* | 1 5 15 35 70 126 210…. | $\frac{x}{(1-x)^5}$ | *-27/1440* | *-27/1440* |
| *6* | *-863/60480* | 1 6 21 56 126 252 462…. | $\frac{x}{(1-x)^6}$ | *-863/60480* | *-863/60480* |
| *7* | *-1375/120960* | 1 7 28 84 210 462 924…. | $\frac{x}{(1-x)^7}$ | *-1375/12096* | *-1375/12096* |

# 4 Ramanujan's Intuitive Method Generalised

Ramanujan's intuitive method to obtain the Ramanujan's sum can be applied to other divergent infinite series $\sum_{n=0}^{\infty} a_n$ that can be represented as $\lim_{x\to 1}\sum_{n=0}^{\infty} a_n x^n$ and generated by algebraic analytic functions $\frac{xF(x)}{(1-x)^k}$ which have a singularity at $x=1^-$ and with *F(x)* being a polynomial. Ramanujan's intuitive method involves using the difference of this function with its twisted or alternating series counterpart $\frac{xF(-x)}{(1+x)^k}$ which has an Abel sum at $x=1^-$. The difference can give the squared counterpart and is represented by the following algebraic equation:

$$\frac{xF(x)}{(1-x)^k}-\frac{xF(-x)}{(1+x)^k}=\frac{x^2F\left(x^2\right)}{\left(1-x^2\right)^k}a_1+x^2H\left(x^2\right) \qquad F(1)\neq 0 \qquad F(-1)\neq 0 \qquad k\in\mathbb{Z},\ \ F(x)\in\mathbb{C}[x] \tag{60}$$

The function $H(x^2)$ must be such that there is a known calculable Ramanujan summation denoted here as $a_2$. The limits as x approaches unity from the left defines the Ramanujan summation:

$$\lim_{x\to1^-}\left[\frac{xF(x)}{(1-x)^k}\right]\overset{\Re}{=}c+\Theta_1 \qquad \lim_{x\to1^-}\left[\frac{x^2F\left(x^2\right)}{\left(1-x^2\right)^k}a_1+x^2H\left(x^2\right)\right]\overset{\Re}{=}a_1c+a_2+\Theta_2 \qquad \lim_{x\to1^-}\left[\frac{xF(-x)}{(1+x)^k}\right]\overset{A}{=}\left[\frac{F(-1)}{2^k}\right]$$

$$\therefore\ \ c+\Theta_1-\left[\frac{F(-1)}{2^k}\right]=a_1c+a_2+\Theta_2 \qquad \therefore\ \ c\overset{\Re}{=}\frac{\left[\frac{F(-1)}{2^k}\right]+a_2}{1-a_1} \tag{61}$$

The variable c is the Ramanujan sum while $\Theta_1$ & $\Theta_2$ are divergent quantities and $a_1$ & $a_2$ are constants. If the divergent quantities are ignored and $a_1$ & $a_2$ are known, then the Ramanujan sum can be calculated. An example is Eqn (10) where $F(x)=1$, $H(x^2)=0,\ a_1=4$ & $a_2=0$.

# 5 Summary

This paper detailed Ramanujan's initial intuitive methodology originally sketched in his notebooks, for assigning the finite value -1/12 to the sum of the divergent infinite series of all positive integers. The methodology was set into an algebraic framework which could be applied to other restricted divergent series. The algebraic framework avoided analytical continuation. The methodology was applied to the infinite sequences of integers raised to a positive integer power and to Figurate Binomial number sequences, including triangular numbers, tetrahedral numbers, and higher-dimensional analogues. The paper applies asymptotic expansions of associated rational generating functions to demonstrate consistency with established results from analytic continuation methods. The results illuminate the deeper mathematical structures underlying Ramanujan's intuitive insights and suggest new avenues for research in divergent series summation.
A key finding was that the regularization of the divergent series involving integers raised to an integer power was shown to be:

$$\lim_{x\to1^-}\left[\frac{xP_{k-1}(x)}{(1-x)^k}-\frac{(k-1)!}{\left(-\ln(x)\right)^k}\right]=\lim_{x\to1^-}\left[\sum_{n=1}^{\infty}n^{k-1}x^n-\frac{(k-1)!}{\left(-\ln(x)\right)^k}\right]=\frac{(-1)^{(k-1)}B_k}{k} \qquad k=1,2,3... \tag{62}$$

The Ramanujan smoothed sums of Figurate Binomial sequences was shown to be intrinsically connected to logarithmic numbers (Gregory coefficients), providing a novel perspective on Ramanujan summation through the lens of classical combinatorial functions.

$$\lim_{x\to1^-}\left[\frac{x}{(1-x)^k}\right]=\lim_{x\to1^-}\left[\sum_{n=1}^{\infty}\binom{n+k-1}{n}x^n\right]\overset{\Re}{=}(-1)^kG_k \qquad k=1,2,3.. \tag{63}$$

By examining the reciprocal series of divergent series using a simple matrix inversion technique, it was shown that the multiplicative inverse or reciprocal of Riemann Zeta function $\sum_{n=1}^{\infty}\frac{1}{n^s}$ with s being positive integers, could be represented by:

$$\sum_{n=1}^{\infty}\frac{1}{n^s}x^{n-1} = \frac{1}{\sum_{n=0}^{\infty}\tilde{G}_n^s x^n} \qquad \tilde{G}_n^s = -\sum_{k=1}^{n}\frac{\tilde{G}_{n-k}^s}{\left(k+1\right)^s} \; n>0 \;\; with\; \tilde{G}_0^s = 1 \tag{64}$$

With coefficients $\tilde{G}_n^s$ being an extended Gregory coefficient system with the denominators being the Hirzebruch numbers raised to the power "s", $\left(h_k\right)^s$.

This work shows that Ramanujan's extraordinary intuition, though initially seeming purely heuristic, embodies deep mathematical substance that continues to reveal new insights under modern analytical scrutiny. As Hardy famously noted, Ramanujan's instincts frequently anticipated results whose rigorous proofs took decades to establish. This paper examined a very simple intuitive idea of Ramanujan and extended it beyond the sum of the positive integers.

# Appendix A– Generating Functions using Linear Algebra

We begin by looking at Euler's alternating series

$$\sum_{n=1}^{\infty}(-1)^{n+1}n \quad = \left(1-2+3-4+5-6+7-8+9-10+\ldots\right) \qquad \mathrm{n}\in\mathbb{N} \tag{65}$$

We list the integer values along a row of a matrix with infinite row length and infinite column length. Here the rows and columns are given a finite length of 9x9 elements. The integer values in each row are staggered by one cell from the previous row, creating an upper triangular matrix. This matrix configuration can easily be inverted without needing any knowledge or interaction from elements at the far end of the row. Only local knowledge is required. The matrix is multiplied by a unit vector and set equal to a vector of constant values denoted by "c".

$$\begin{bmatrix} 1 & -2 & 3 & -4 & 5 & -6 & 7 & -8 & 9\ldots \\ 0 & 1 & -2 & 3 & -4 & 5 & -6 & 7 & -8\ldots \\ 0 & 0 & 1 & -2 & 3 & -4 & 5 & -6 & 7\ldots \\ 0 & 0 & 0 & 1 & -2 & 3 & -4 & 5 & -6\ldots \\ 0 & 0 & 0 & 0 & 1 & -2 & 3 & -4 & 5\ldots \\ 0 & 0 & 0 & 0 & 0 & 1 & -2 & 3 & -4\ldots \\ 0 & 0 & 0 & 0 & 0 & 0 & 1 & -2 & 3\ldots \\ 0 & 0 & 0 & 0 & 0 & 0 & 0 & 1 & -2\ldots \\ 0 & 0 & 0 & 0 & 0 & 0 & 0 & 0 & 1 \end{bmatrix} \begin{bmatrix} 1 \\ 1 \\ 1 \\ 1 \\ 1 \\ 1 \\ 1 \\ 1 \\ 1 \end{bmatrix} = \begin{bmatrix} c \\ c \\ c \\ c \\ c \\ c \\ c \\ c \\ c \end{bmatrix} \tag{66}$$

The inverse is:

$$\begin{bmatrix} 1 & 2 & 1 & 0 & 0 & 0 & 0 & 0 & 0\ldots \\ 0 & 1 & 2 & 1 & 0 & 0 & 0 & 0 & 0\ldots \\ 0 & 0 & 1 & 2 & 1 & 0 & 0 & 0 & 0\ldots \\ 0 & 0 & 0 & 1 & 2 & 1 & 0 & 0 & 0\ldots \\ 0 & 0 & 0 & 0 & 1 & 2 & 1 & 0 & 0\ldots \\ 0 & 0 & 0 & 0 & 0 & 1 & 2 & 1 & 0\ldots \\ 0 & 0 & 0 & 0 & 0 & 0 & 1 & 2 & 1\ldots \\ 0 & 0 & 0 & 0 & 0 & 0 & 0 & 1 & 2\ldots \\ 0 & 0 & 0 & 0 & 0 & 0 & 0 & 0 & 1 \end{bmatrix} \begin{bmatrix} c \\ c \\ c \\ c \\ c \\ c \\ c \\ c \\ c \end{bmatrix} = \begin{bmatrix} 1 \\ 1 \\ 1 \\ 1 \\ 1 \\ 1 \\ 1 \\ 1 \\ 1 \end{bmatrix} \tag{67}$$

Assuming the series can be written as the coefficients of a polynomial for example:

$$\sum_{n=1}^{\infty}(-1)^{n+1}nx^{n-1} \quad = \left(1-2x+3x^2-4x^3+5x^4-6x^5+7x^6-8x^7+9x^8-10x^9+\ldots\right) \qquad \mathrm{n}\in\mathbb{N} \tag{68}$$

then Eqn (67) gives us information to construct the multiplicative inverse generating function, that is:

$$\left(1-2x+3x^2-4x^3+5x^4-6x^5+7x^6-8x^7+\ldots\right) = \left(\frac{1}{1+2x+x^2} = \frac{1}{\left(1+x\right)^2}\right) \tag{69}$$

The second example is the alternating unit series:

$$(1-1+1-1+1-1+1-1+1-1+...) \tag{70}$$

Again, writing in 9x9 matrix form, we have:

$$\begin{bmatrix} 1 & -1 & 1 & -1 & 1 & -1 & 1 & -1 & 1... \\ 0 & 1 & -1 & 1 & -1 & 1 & -1 & 1 & -1... \\ 0 & 0 & 1 & -1 & 1 & -1 & 1 & -1 & 1... \\ 0 & 0 & 0 & 1 & -1 & 1 & -1 & 1 & -1... \\ 0 & 0 & 0 & 0 & 1 & -1 & 1 & -1 & 1... \\ 0 & 0 & 0 & 0 & 0 & 1 & -1 & 1 & -1... \\ 0 & 0 & 0 & 0 & 0 & 0 & 1 & -1 & 1... \\ 0 & 0 & 0 & 0 & 0 & 0 & 0 & 1 & -1... \\ 0 & 0 & 0 & 0 & 0 & 0 & 0 & 0 & 1... \end{bmatrix} \begin{bmatrix} 1 \\ 1 \\ 1 \\ 1 \\ 1 \\ 1 \\ 1 \\ 1 \\ 1 \end{bmatrix} = \begin{bmatrix} c \\ c \\ c \\ c \\ c \\ c \\ c \\ c \\ c \end{bmatrix} \tag{71}$$

The inverse is:

$$\begin{bmatrix} 1 & 1 & 0 & 0 & 0 & 0 & 0 & 0 & 0... \\ 0 & 1 & 1 & 0 & 0 & 0 & 0 & 0 & 0... \\ 0 & 0 & 1 & 1 & 0 & 0 & 0 & 0 & 0... \\ 0 & 0 & 0 & 1 & 1 & 0 & 0 & 0 & 0... \\ 0 & 0 & 0 & 0 & 1 & 1 & 0 & 0 & 0... \\ 0 & 0 & 0 & 0 & 0 & 1 & 1 & 0 & 0... \\ 0 & 0 & 0 & 0 & 0 & 0 & 1 & 1 & 0... \\ 0 & 0 & 0 & 0 & 0 & 0 & 0 & 1 & 1... \\ 0 & 0 & 0 & 0 & 0 & 0 & 0 & 0 & 1... \end{bmatrix} \begin{bmatrix} c \\ c \\ c \\ c \\ c \\ c \\ c \\ c \\ c \end{bmatrix} = \begin{bmatrix} 1 \\ 1 \\ 1 \\ 1 \\ 1 \\ 1 \\ 1 \\ 1 \\ 1 \end{bmatrix} \tag{72}$$

Representing the series as a polynomial, thus:

$$\sum_{n=1}^{\infty} (-1)^{n-1} x^{n-1} = \left(1-x+x^2-x^3+x^4-x^5+x^6-x^7+x^8-x^9+...\right) \quad \text{n} \in \mathbb{N} \tag{73}$$

From Eqn (72) we conclude that:

$$\left(1-x+x^2-x^3+x^4-x^5+x^6-x^7+x^8-x^9+...\right) = \left(\frac{1}{1+x}\right) \tag{74}$$

The third example is the Grandi's series:

$$(1+1+1+1+1+1+1+1+1+1+...) \tag{75}$$

Writing in matrix form as:

$$\begin{bmatrix} 1 & 1 & 1 & 1 & 1 & 1 & 1 & 1 & 1\ldots \\ 0 & 1 & 1 & 1 & 1 & 1 & 1 & 1 & 1\ldots \\ 0 & 0 & 1 & 1 & 1 & 1 & 1 & 1 & 1\ldots \\ 0 & 0 & 0 & 1 & 1 & 1 & 1 & 1 & 1\ldots \\ 0 & 0 & 0 & 0 & 1 & 1 & 1 & 1 & 1\ldots \\ 0 & 0 & 0 & 0 & 0 & 1 & 1 & 1 & 1\ldots \\ 0 & 0 & 0 & 0 & 0 & 0 & 1 & 1 & 1\ldots \\ 0 & 0 & 0 & 0 & 0 & 0 & 0 & 1 & 1\ldots \\ 0 & 0 & 0 & 0 & 0 & 0 & 0 & 0 & 1\ldots \end{bmatrix} \begin{bmatrix} 1 \\ 1 \\ 1 \\ 1 \\ 1 \\ 1 \\ 1 \\ 1 \\ 1 \end{bmatrix} = \begin{bmatrix} c \\ c \\ c \\ c \\ c \\ c \\ c \\ c \\ c \end{bmatrix} \tag{76}$$

The inverse matrix is:

$$\begin{bmatrix} 1 & -1 & 0 & 0 & 0 & 0 & 0 & 0 & 0\ldots \\ 0 & 1 & -1 & 0 & 0 & 0 & 0 & 0 & 0\ldots \\ 0 & 0 & 1 & -1 & 0 & 0 & 0 & 0 & 0\ldots \\ 0 & 0 & 0 & 1 & -1 & 0 & 0 & 0 & 0\ldots \\ 0 & 0 & 0 & 0 & 1 & -1 & 0 & 0 & 0\ldots \\ 0 & 0 & 0 & 0 & 0 & 1 & -1 & 0 & 0\ldots \\ 0 & 0 & 0 & 0 & 0 & 0 & 1 & -1 & 0\ldots \\ 0 & 0 & 0 & 0 & 0 & 0 & 0 & 1 & -1\ldots \\ 0 & 0 & 0 & 0 & 0 & 0 & 0 & 0 & 1\ldots \end{bmatrix} \begin{bmatrix} c \\ c \\ c \\ c \\ c \\ c \\ c \\ c \\ c \end{bmatrix} = \begin{bmatrix} 1 \\ 1 \\ 1 \\ 1 \\ 1 \\ 1 \\ 1 \\ 1 \\ 1 \end{bmatrix} \tag{77}$$

Converting to polynomial form, we have:

$$\sum_{n=1}^{\infty} x^{n-1} \quad = \left(1+x+x^2+x^3+x^4+x^5+x^6+x^7+x^8+x^9+\ldots\right) \qquad \mathrm{n} \in \mathbb{N} \tag{78}$$

From this we conclude that:

$$\left(1+x+x^2+x^3+x^4+x^5+x^6+x^7+x^8+x^9+\ldots\right)=\left(\frac{1}{1-x}\right) \tag{79}$$

The last example relates to the Harmonic series $\sum_{n=1}^{\infty}\frac{1}{n}$ written in matrix form as:

$$\begin{bmatrix} 1 & \frac{1}{2} & \frac{1}{3} & \frac{1}{4} & \frac{1}{5} & \frac{1}{6} & \frac{1}{7} & \frac{1}{8} & \frac{1}{9}\ldots \\ 0 & 1 & \frac{1}{2} & \frac{1}{3} & \frac{1}{4} & \frac{1}{5} & \frac{1}{6} & \frac{1}{7} & \frac{1}{8}\ldots \\ 0 & 0 & 1 & \frac{1}{2} & \frac{1}{3} & \frac{1}{4} & \frac{1}{5} & \frac{1}{6} & \frac{1}{7}\ldots \\ 0 & 0 & 0 & 1 & \frac{1}{2} & \frac{1}{3} & \frac{1}{4} & \frac{1}{5} & \frac{1}{6}\ldots \\ 0 & 0 & 0 & 0 & 1 & \frac{1}{2} & \frac{1}{3} & \frac{1}{4} & \frac{1}{5}\ldots \\ 0 & 0 & 0 & 0 & 0 & 1 & \frac{1}{2} & \frac{1}{3} & \frac{1}{4}\ldots \\ 0 & 0 & 0 & 0 & 0 & 0 & 1 & \frac{1}{2} & \frac{1}{3}\ldots \\ 0 & 0 & 0 & 0 & 0 & 0 & 0 & 1 & \frac{1}{2}\ldots \\ 0 & 0 & 0 & 0 & 0 & 0 & 0 & 0 & 1 \end{bmatrix} \begin{bmatrix} 1 \\ 1 \\ 1 \\ 1 \\ 1 \\ 1 \\ 1 \\ 1 \\ 1 \end{bmatrix} = \begin{bmatrix} c \\ c \\ c \\ c \\ c \\ c \\ c \\ c \\ c \end{bmatrix} \tag{80}$$

The inverse is given by:

$$\begin{bmatrix} 1 & \frac{-1}{2} & \frac{-1}{12} & \frac{-1}{24} & \frac{-19}{720} & \frac{-3}{160} & \frac{-863}{60480} & \frac{-275}{24192} & \ldots \\ 0 & 1 & \frac{-1}{2} & \frac{-1}{12} & \frac{-1}{24} & \frac{-19}{720} & \frac{-3}{160} & \frac{-863}{60480} & \ldots \\ 0 & 0 & 1 & \frac{-1}{2} & \frac{-1}{12} & \frac{-1}{24} & \frac{-19}{720} & \frac{-3}{160} & \ldots \\ 0 & 0 & 0 & 1 & \frac{-1}{2} & \frac{-1}{12} & \frac{-1}{24} & \frac{-19}{720} & \ldots \\ 0 & 0 & 0 & 0 & 1 & \frac{-1}{2} & \frac{-1}{12} & \frac{-1}{24} & \ldots \\ 0 & 0 & 0 & 0 & 0 & 1 & \frac{-1}{2} & \frac{-1}{12} & \ldots \\ 0 & 0 & 0 & 0 & 0 & 0 & 1 & \frac{-1}{2} & \ldots \\ 0 & 0 & 0 & 0 & 0 & 0 & 0 & 1 & \ldots \\ 0 & 0 & 0 & 0 & 0 & 0 & 0 & 0 & 1 \end{bmatrix} \begin{bmatrix} c \\ c \\ c \\ c \\ c \\ c \\ c \\ c \\ c \end{bmatrix} = \begin{bmatrix} 1 \\ 1 \\ 1 \\ 1 \\ 1 \\ 1 \\ 1 \\ 1 \\ 1 \end{bmatrix} \tag{81}$$

The coefficients in the rows of the inverse matrix Eqn (81) are related to the Gregory coefficients. Assuming the series can be written as the coefficients of a polynomial, we have:

$$\begin{aligned} \sum_{n=1}^{\infty} \frac{1}{n} x^{n-1} &= \left(1 + \frac{1}{2}x + \frac{1}{3}x^2 + \frac{1}{4}x^3 + \frac{1}{5}x^4 + \frac{1}{6}x^5 + \frac{1}{7}x^6 + \frac{1}{8}x^7 + \frac{1}{9}x^8 + \frac{1}{10}x^9 + \ldots\right) \\ = \frac{1}{\sum_{n=0}^{\infty} (-1)^n G_n x^n} &= \frac{1}{\left(1 - \frac{1}{2}x - \frac{1}{12}x^2 - \frac{1}{24}x^3 - \frac{19}{720}x^4 - \frac{27}{1440}x^5 - \frac{863}{60480}x^6 - \frac{1375}{120960}x^7 - \ldots\right)} \end{aligned} \tag{82}$$

As detailed earlier in Eqn (44), the Gregory coefficients are related to the Harmonic series. The recursive formula for the Gregory coefficients, Eqn (45) is rewritten here:

$$\bar{G}_n = -\sum_{k=1}^{n} \frac{\bar{G}_{n-k}}{(k+1)} \quad n>0 \;\; with \;\; \bar{G}_0 = 1 \qquad G_0 = 1 \quad G_n = (-1)^n \bar{G}_n \quad n = 1,2,3... \tag{83}$$

As stated earlier, the denominators of the Gregory coefficients are the Hirzebruch numbers [17, 18]:

$$h_k = 1,\ 2,\ 12,\ 24,\ 720,\ 1440,\ 60480,.... = \prod_{prime\ p \le k+1} p^{\left\lfloor \frac{k}{p-1} \right\rfloor} \tag{84}$$

For the Riemann Zeta function $\sum_{n=1}^{\infty} \frac{1}{n^s}$, the same process can be applied to obtain the multiplicative inverse or reciprocal. The following is obtained:

$$\sum_{n=1}^{\infty} \frac{1}{n^s} x^{n-1} = \frac{1}{\sum_{n=0}^{\infty} \tilde{G}_n^s x^n} \qquad \tilde{G}_n^s = -\sum_{k=1}^{n} \frac{\tilde{G}_{n-k}^s}{(k+1)^s} \;\; n>0 \;\; with \;\; \tilde{G}_0^s = 1 \tag{85}$$

With coefficients $\tilde{G}_n^s$ being an extended Gregory coefficient system with the denominators being the Hirzebruch numbers raised to the power "s", $(h_k)^s$. Examples of the extended Greogory coefficients are given in Table 5.

Table 5 Examples of Extended Gregory coefficients for s=1 to 5

| $n$ | 1 | 2 | 3 | 4 | 5 | 6 | 7 |
|---|---|---|---|---|---|---|---|
| $\tilde{G}_n^1$ | 1 | -1/2 | -1/12 | -1/24 | -19/720 | -27/1440 | -863/60480 |
| $\tilde{G}_n^2$ | 1 | $-1/2^2$ | $-7/12^2$ | $-13/24^2$ | $-6911/720^2$ | $-18453/1440^2$ | $-23419855/60480^2$ |
| $\tilde{G}_n^3$ | 1 | $-1/2^3$ | $-37/12^3$ | $-115/24^3$ | $-1572859/720^3$ | $-7346157/1440^3$ | $-347737791311/60480^3$ |
| $\tilde{G}_n^4$ | 1 | $-1/2^4$ | $-175/12^4$ | $-865/24^4$ | $-292581071/720^4$ | $-2315047233/1440^4$ | $-3980321414135551/60480^4$ |
| $\tilde{G}_n^5$ | 1 | $-1/2^5$ | $-781/12^5$ | $-5971/24^5$ | $-48979036099/720^5$ | $-215760166559*3/1440^5$ | $-40003092470353818383/60480^5$ |

The first six values for the extended Gregory coefficients are given below in symbolic form.

$$\tilde{G}_0^s = 1 \quad \tilde{G}_1^s = -\frac{1}{2^s} \quad \tilde{G}_2^s = \frac{3^s - \left(2^s\right)^2}{\left(2^2 3\right)^s} = \frac{1}{\left(2^2\right)^s} - \frac{1}{3^s} \quad \tilde{G}_3^s = \frac{\left(2\left(2^s\right)^2 - 3^s\right)4^s - \left(2^s\right)^3 3^s}{\left(2^3 3.4\right)^s} = \frac{2}{2^s 3^s} - \frac{1}{\left(2^3\right)^s} - \frac{1}{4^s}$$

$$\tilde{G}_4^s = \frac{-\left(2^s\right)^4\left(\left(3^s\right)^2 - 5^s\right)4^s + 2\left(2^s\right)^3\left(3^s\right)^2 5^s - 3\left(2^s\right)^2 3^s 4^s 5^s + \left(3^s\right)^2 4^s 5^s}{\left(2^4 3^2 4.5\right)^s} = \frac{1}{\left(2^4\right)^s} - \frac{3}{\left(2^2\right)^s 3^s} + \frac{2}{2^s 4^s} + \frac{1}{3^s} - \frac{1}{5^s} \tag{86}$$

$$\tilde{G}_5^s = \frac{-\left(2^s\right)^5\left(3^s 4^s - 2.6^s\right)3^s 5^s + \left(2^s\right)^4\left(2\left(3^s\right)^2 - 3.5^s\right)4^s 6^s - 3\left(2^s\right)^3\left(3^s\right)^2 5^s 6^s + 4\left(2^s\right)^2 3^s 4^s 5^s 6^s - \left(3^s\right)^2 4^s 5^s 6^s}{\left(2^5 3^2 4.5.6\right)^s}$$

$$= -\frac{1}{\left(2^5\right)^s} + \frac{4}{\left(2^3\right)^s 3^s} - \frac{3}{\left(2^2\right)^s 4^s} - \frac{3}{2^s\left(3^2\right)^s} + \frac{2}{3^s 4^s} + \frac{2}{2^s 5^s} - \frac{1}{6^s}$$

# Appendix B Higher Order Derivatives of Eulerian Polynomials evaluated at x=1

The first four derivatives of the first seven Eulerian Polynomials are listed below:

| $u$ | $P_u^{'}(x)$ | $P_u^{''}(x)$ | $P_u^{'''}(x)$ | $P_u^{''''}(x)$ |
|---|---|---|---|---|
| 0 | 0 | 0 | 0 | 0 |
| 1 | 0 | 0 | 0 | 0 |
| 2 | 1 | 0 | 0 | 0 |
| 3 | $4+2x$ | 2 | 0 | 0 |
| 4 | $11+22x+3x^2$ | $22+6x$ | 6 | 0 |
| 5 | $26+132x+78x^2+4x^3$ | $132+156x+12x^2$ | $156+24x$ | 24 |
| 6 | $57+604x+906x^2+228x^3+5x^4$ | $604+1812x+684x^2+20x^3$ | $1812+1368x+60x^2$ | $1368+120x$ |

(87)

The following identity for the $m^{th}$ derivative of $u^{th}$ Eulerian polynomials evaluated at *x=1* are:

$$\left(\frac{d^m P_u(x)}{dx^m}\right)_{x=1} = (-1)^m u!(u-m)m!G_m^{(u)} = u!B_m^{(m-u)} \quad u > m \qquad m = 0,1,2,3,4...$$

$$\left(\frac{d^m P_u(x)}{dx^m}\right) = 0 \quad u \le m \tag{88}$$

Where $G_m^{(u)}$ are the $m^{th}$ Gregory polynomials in the variable *u* defined in Appendix C. Equation (97) has been used to write the derivatives in terms of the Bernoulli polynomial. Expanding the above for the first five derivatives evaluated at *x=1*, we have:

$$
\begin{aligned}
&m=0 \quad \left(P_u(x)\right)_{x=1} &&= u! \\
&m=1 \quad \left(\frac{dP_u(x)}{dx}\right)_{x=1} &&= u!(u-1)\frac{1}{2} \\
&m=2 \quad \left(\frac{d^2P_u(x)}{dx^2}\right)_{x=1} &&= u!(u-2)\frac{(3u-5)}{12} \\
&m=3 \quad \left(\frac{d^3P_u(x)}{dx^3}\right)_{x=1} &&= u!\frac{(u-2)(u-3)^2}{8} \\
&m=4 \quad \left(\frac{d^4P_u(x)}{dx^4}\right)_{x=1} &&= u!(u-4)\frac{\left(15u^3-150u^2+485u-502\right)}{240} \\
&m=5 \quad \left(\frac{d^5P_u(x)}{dx^5}\right)_{x=1} &&= u!\frac{(u-5)^2(u-4)\left(3u^2-23u+38\right)}{96}
\end{aligned}
\tag{89}
$$

# Appendix C Gregory Polynomials

In this paper we define the $m^{th}$ Gregory polynomial in the variable *u,* as $G_m^{(u)}$. The polynomial can be extracted from the following series expansion:

$$
\begin{aligned}
\frac{\left(\dfrac{-x}{\ln(1-x)}\right)^u}{u} &= \frac{1}{u}-\frac{1}{2}x+\frac{(3u-5)}{24}x^2-\frac{(u-2)(u-3)}{48}x^3+\frac{\left(15u^3-150z^2++485u-502\right)}{5760}x^4 \\
&-\frac{(u-5)(u-4)\left(3u^2-23u+38\right)}{11520}x^5+\frac{\left(63u^5-1575u^4+15435u^3-73801u^2+171150u-152696\right)}{2903040}x^6-..... \\
&=\sum_{m=0}^{\infty} G_m^{(u)}x^m
\end{aligned}
\tag{90}
$$

The above generating function is discussed in OEIS A075264. Note from Eqn (44) that $\left(\dfrac{x}{-\ln(1-x)}\right)$ is the generating function for the Gregory coefficients. Therefore, the polynomials $G_m^{(u)}$ could be referred to as "Gregory polynomials" in the variable *u*. The polynomials can be obtained using the recursive formula:

$$
\begin{aligned}
&\bar{G}_0^{(u)}=1 \quad \bar{G}_1^{(u)}=\frac{-u}{2} \quad \bar{G}_{i+1}^{(u)}=-\frac{u}{i+2}-\frac{1}{i+1}\left(\sum_{n=1}^{i}\frac{u(i-n+1)+n}{(i-n+2)}\bar{G}_n^{(u)}\right) \quad i=1,2,3... \\
&G_m^{(u)}=\frac{\bar{G}_m^{(u)}}{u} \quad m=0,1,2,3....
\end{aligned}
\tag{91}
$$

The Gregory coefficients emerge from $G_m=(-1)^m G_m^{(1)}$ & $\frac{1}{m}=-G_m^{(-1)}$. The Gregory polynomials also satisfy the following identity:

$$G_m^{(u+1)} = \frac{\left[(u-m)G_m^{(u)} - (u-m+1)G_{m-1}^{(u)}\right]}{(u+1)} \tag{92}$$

Substituting *m=k & u=k-1* into the above equation gives the identity:

$$k\,G_k^{(k)} + G_k^{(k-1)} = 0 \tag{93}$$

The first six Gregory polynomials as functions of the variable "u" are written below.

$$\begin{aligned}
m &= 0 \quad G_0^{(u)} = \frac{1}{u} \\
m &= 1 \quad G_1^{(u)} = -\frac{1}{2} \\
m &= 2 \quad G_2^{(u)} = \frac{3u-5}{24} \\
m &= 3 \quad G_3^{(u)} = -\frac{(u-2)(u-3)}{48} \\
m &= 4 \quad G_4^{(u)} = \frac{15u^3 - 150u^2 + 485u - 502}{5760} \\
m &= 5 \quad G_5^{(u)} = -\frac{(u-5)(u-4)\left(3u^2 - 23u + 38\right)}{11520}
\end{aligned} \tag{94}$$

The Gregory polynomials can be written using either of the following generating functions:

$$\left(\frac{x}{-\ln(1-x)}\right)^u = u\sum_{m=0}^{\infty} G_m^{(u)} x^m \qquad \left(\frac{\ln(1+x)}{x}\right)^u = u\sum_{m=0}^{\infty} (-1)^m G_m^{(-u)} x^m \tag{95}$$

The following equation is detailed in references [21, 22].

$$\left(\frac{\ln(1+x)}{x}\right)^{-u} = -u\sum_{m=0}^{\infty} \frac{B_m^{(m-u)}}{(m-u)} \frac{x^m}{m!} \qquad |x| < 1 \tag{96}$$

This equation can be used to derive a relationship between the Bernoulli numbers and Gregory polynomials. Using the above and Eqn (93) and (95), we can write:

$$\left(\frac{\ln(1+x)}{x}\right)^{-u} = -u\sum_{m=0}^{\infty} \frac{B_m^{(m-u)}}{(m-u)} \frac{x^m}{m!} = -u\sum_{m=0}^{\infty} (-1)^m G_m^{(u)} x^m \qquad \therefore \quad B_m^{(m-u)} = (-1)^{m+1}(m-u)\,m!\,G_m^{(u)}$$

$$let \quad m = k = 1+u \qquad \therefore \quad \frac{B_k}{k!} = (-1)^{k+1} G_k^{(k-1)} = (-1)^k\,k\,G_k^{(k)} \tag{97}$$

Also, we can derive the following expression with *u=1,* in Eqn (97) (see [22] Eqn 2.12).

$$\frac{B_m^{(m-1)}}{(m-1)\,m!} = (-1)^{m+1} G_m^{(1)} = -G_m \tag{98}$$